\documentclass[11pt,reqno]{amsproc}

\title[]{Lagrangian-Eulerian Methods for Uniqueness in Hydrodynamic Systems}

\author{Peter Constantin}
\address{Department of Mathematics, Princeton University, Princeton, NJ 08544}
\email{const@math.princeton.edu}

\usepackage[margin=1in]{geometry}
\usepackage{amsmath, amsthm, amssymb}
\usepackage{times}

\usepackage{color}
\usepackage{hyperref}

\newcommand{\pa}{\partial}
\newcommand{\la}{\label}
\newcommand{\fr}{\frac}
\newcommand{\na}{\nabla}
\newcommand{\be}{\begin{equation}}
\newcommand{\ee}{\end{equation}}
\newcommand{\ba}{\begin{array}{l}}
\newcommand{\ea}{\end{array}}
\newcommand{\Rr}{{\mathbb R}}

\newcommand{\beg}{\begin}

\renewcommand{\div}{{\mbox{div}\,}}
\date{today}
\begin{document}
\begin{abstract}
We present a Lagrangian-Eulerian strategy for proving uniqueness and local existence of solutions of limited smoothness for a class of incompressible hydrodynamic models including Oldroyd-B type complex fluid models and zero magnetic resistivity magneto-hydrodynamics equations.
\end{abstract}

\maketitle
\section{Introduction}
Many physical models consist of equations for fluids coupled with equations for other fields. Primary examples occur in descriptions of complex fluids in which a solvent interacts with particles, and  in magneto-hydrodynamics, in which a fluid interacts with a magnetic field. One of the simplest complex fluids models, an Oldroyd-B model, reduces to a time independent Stokes system
\[
-\Delta u + \na p = \div \sigma, \quad \div u = 0
\]
coupled with an evolution equation for the symmetric added stress matrix $\sigma$,
\[
\pa_t \sigma + u\cdot\na\sigma = (\na u)\sigma + \sigma (\na u)^T - \sigma + (\na u) + (\na u)^T
\]
Clearly, from the Stokes equation with appropriate boundary conditions (for instance decay in the whole space) it follows that the velocity gradient is of the same order of magnitude as the added stress, $\na u\sim \sigma$. This makes the evolution equation for $\sigma$ potentially capable of producing finite time blow up. The formation of finite time singularities in this system is an outstanding open problem. While the balance $\sigma\sim\na u$ is potentially dangerous for large data and long time, it also indicates clearly that if $\na u$ is controlled then $\sigma$ is controlled as well. In particular, the short time existence of solutions can be obtained in a class of velocities that is close to the Lipschitz class. The fact that singular integral operators are not bounded in $L^{\infty}$ requires the use of slightly smaller spaces, and $\sigma\in C^{\alpha}(\Rr^d)\cap L^p(\Rr^d)$ and
correspondingly $u\in C^{1+\alpha}(\Rr^d)\cap W^{1,p}(\Rr^d)$ are spaces in which the problem admits short time existence. Here both $\alpha\in (0,1)$ and $p\in (1,\infty)$ are arbitrary. It is natural then to ask about uniqueness of solutions in the same spaces. Taking the difference $\sigma$ between two solutions $\sigma_1$ and $\sigma_2$ leads to an equation 
\[
%\ba
\pa_t \sigma + {\bar u}\cdot\na\sigma + u\cdot\na {\bar{\sigma}} =
(\na \bar{u})\sigma + (\na u)\bar{\sigma} + \bar\sigma(\na u)^T + \sigma (\na\bar u)^T -\sigma + (\na u) + (\na u)^T
%\ea
\]
where $u = u_1-u_2$ is the difference of the corresponding velocities, and $\bar u = \fr{1}{2}(u_1+u_2)$ and $\bar\sigma = \fr{1}{2}(\sigma_1+ \sigma_2)$ are the arithmetic averages of velocities and of stresses. The right hand side is well-behaved in $C^{\alpha}$. The term $u\cdot\na\bar\sigma$ is not defined for $\bar\sigma \in C^{\alpha}$.
This makes an Eulerian approach to a uniqueness proof difficult in this class of solutions. Uniqueness with this low regularity was proved in (\cite{csun}), using however a Lagrangian approach. The main reason why Lagrangian variables are better behaved than Eulerian ones is that in Lagrangian variables the velocity $v= u\circ X$ is obtained from the Lagrangian added stress $\tau = \sigma\circ X$ by an expression
\[
v = {\mathbb U}(\tau\circ X^{-1})\circ X
\]
where ${\mathbb U}$ is the linear operator that produces the solution of the steady Stokes equation from the added stresses, and $X$ is the Lagrangian path, which is a time-dependent diffeomorphism. The Gateaux derivative (variational derivative or first variation, in the language of mechanics) of the map $X\mapsto v$ is a commutator, and it is better behaved than each of its terms. On the other hand, $\tau$ obeys an ODE in Lagrangian varaibles, so it is easily controlled for short time by $g= (\na u)\circ X$. 

The present paper expands this approach to time-dependent relationships between $u$ and $\sigma$ and we prove uniqueness and local existence in large spaces $C^{\alpha}\cap L^p$ for a class of hydrodynamic models including complex fluids of Oldroyd-B type, and ideal magneto-hydrodynamics. Local existence of very smooth solutions of such systems is classical (\cite{saut}). We provide in this paper a correct statement and a complete proof of a lemma (Lemma \ref{steadycomm}) which was used in (\cite{csun}) quoting (\cite{mb}). A Lagrangian approach for Oldroyd-B smooth flow was advocated in (\cite{fanghuachun}). Adding inertia, i.e. coupling with Navier-Stokes instead of Stokes, requires a modified treatment. We prove general existence and uniqueness theorems in Lagarngian variables. They apply in particular to the ideal MHD equations 
\[
\left\{
\ba
\pa_t u + u\cdot\na -\nu\Delta u + \na p = b\cdot\na b,\\
\div u = 0,\\
\pa_t b + u\cdot\na b = b\cdot\na u
\ea
\right.
\]
in ${\Rr^d}$, and to nonlinear Oldroyd-B like systems
\[
\left\{
\ba
\pa_t u + u\cdot\na -\nu\Delta u + \na p = \div\sigma,\\
\div u = 0,\\
\pa_t \sigma + u\cdot\na \sigma = F(\na u, \sigma)
\ea
\right.
\]
in $\Rr^d$, for quite general smooth $F$. The results prove local existence and Lipschitz dependence on initial data in Lagrangian coordinates, in classes of H\"{o}lder continuous
magnetic field $b$, added stress $\sigma$, and velocity gradient $\na u$.

The paper is organized as follows. In the next section we describe systems in which added stresses are coupled to time dependent Stokes equations. The third section is devoted to statements and proofs about the linear operators and the commutators involved. The fourth section presents the proof of Theorems \ref{laglip} and \ref{lagexist}, for Stokes-based systems, which state that in Lagrangian variables the nonlinearity is locally Lipschitz in path spaces. The solutions exist locally, and depend in a locally Lipschitz continuous manner in these certain path spaces on initial data. In particular, the solutions are unique. Although simpler, the time-dependent Stokes-based systems provide the principal challenges. Once the setting is clarified, the Navier-Stokes-based systems are treated in these path spaces in a perturbative manner. The difference is that the Stokes-based systems close at the level of $(X,\tau)$, whereas the Navier-Stokes-based systems require a treatment with $(X,\tau, \pa_t X)$ simultaneously. The fifth section describes the changes needed in order to adapt the proof for the case of inertia. The main results, Theorems \ref{contrns} and \ref{exns} state local existence of solutions for Navier-Stokes-based systems in Lagrangian variables and Lipschitz dependence of initial data in a path space
\[ 
(X,\tau, \pa_t X)\in C^{\beta}(0,T; C^{1+\alpha, p})\times C^{\beta}(0,T; C^{\alpha, p}(\Rr^d))\times L^{\infty}(0,T; C^{1+\alpha,p}).
\]

\section{Time-dependent Stokes flow coupled with added stress}
We consider the time dependent forced Stokes equations
\be
\pa_t u -\Delta u + \na p = \div\sigma, \quad \na\cdot u = 0
\la{stokes}
\ee
coupled with an equation
\be
\pa_t \sigma + u\cdot\na\sigma = F((\na u), \sigma )
\la{sigeq}
\ee
We assume that $F$ is a smooth function. We also assume that $F$ has polynomial  growth at infinity 
\[
|F(g,\tau)|\le C(|g| +|\tau|)^k.
\]
for some $k\in \mathbb N$. This is only to guarantee some decay at infinity, because we work in $L^p\cap L^{\infty}$ spaces in $\Rr^d$.
The typical examples include
\[
F((\na u), \sigma) = (\na u)\sigma + \sigma (\na u)^T -\sigma + (\na u) + (\na u)^T
\]
In the case of MHD, $\sigma = b$ is a vector, the right-hand side of the Stokes equation is given by $\div(b\otimes b)$ and 
\[
F = b\cdot\na u.
\]
The divergence operation is $\div \sigma = \na\cdot\sigma$.
We work in ${\Rr^d}$ and require the velocities $u$ and the stresses $\sigma$ to vanish at infinity. The solution map for the Stokes equation is
\be
u(x,t) = {\mathbb L}(u_0)(x,t)+ {\mathbb{U}}(\sigma)(x,t)
\la{u}
\ee
where
\be
{\mathbb{U}}(\sigma) = \int_0^te^{(t-s)\Delta}{\mathbb H}\div\sigma (s)ds 
\la{bigu}
\ee 
with
\be
{\mathbb H} = {\mathbb I} + R\otimes R,
\la{h}
\ee
and $R= (-\Delta)^{-\fr{1}{2}}\na$ the Riesz transforms, and with
\be
{\mathbb L}(u_0)(t) =  e^{t\Delta}u_0
\la{bigl}
\ee
Note that ${\mathbb U}(\sigma)$ is divergence-free at vanishes at $t=0$, and
${\mathbb L}(u_0)$ is divergence-free if $u_0$ is. Also,
\be
\na u = e^{t\Delta}\na u_0 + {\mathbb G}(\sigma) = {\mathbb L}(\na u_0) + {\mathbb G}(\sigma)
\la{nau}
\ee
with
\be
{\mathbb G}(\sigma) = \int_0^te^{(t-s)\Delta}\na{\mathbb H}\div\sigma (s)ds.
\la{bigg}
\ee
The Lagrangian description is as follows. The Lagrangian paths $X$ solve
\be
\fr{dX}{dt} = {\mathbb L}(u_0)\circ X + {\mathbb U}( \tau\circ X^{-1})\circ X
\la{xeq}
\ee
with
\be
\tau = \sigma\circ X
\la{tau}
\ee
and initial data $X(a,0)=a$. Let
\be
g(a,t) = (\na u)(X(a,t),t) = {\mathbb L}(\na u_0)\circ X + {\mathbb G}(\tau\circ X^{-1})\circ X.
\la {g}
\ee
The equation for $\tau$ is the ODE
\be
\fr{d\tau}{dt} = F(g,\tau).
\la{tauode}
\ee
The Eulerian variables are $u$ and $\sigma$. The Lagrangian variables are $X$ and $\tau$. In Lagrangian variables, the system is
\be
\left\{
\ba
X(a,t) = a + \int_0^t{\mathcal U}(X(s),\tau(s))ds,\\
\tau (a,t) = \sigma_0(a) + \int_0^t{\mathcal T}(X(s), \tau(s))ds.
\ea
\right.
\la{lagfix}
\ee
where the Lagrangian nonlinearities ${\mathcal U}(X,\tau)$, ${\mathcal T}(X,\tau)$ are
\be
\left\{
\ba
\mathcal U = {\mathbb L}(u_0)\circ X + {\mathbb U}(\tau\circ X^{-1})\circ X,\\
\mathcal T = F\left(\left({\mathbb L}(\na u_0)\circ X + {\mathbb G}(\tau\circ X^{-1})\circ X\right), \tau \right)
\ea
\right .
\la{lagsys}
\ee
We consider a differentiable one-parameter family of paths $X_{\epsilon}, \tau_{\epsilon}$, with Eulerian form $\sigma_{\epsilon} = \tau_{\epsilon}\circ X_{\epsilon}^{-1}$, initial data $u_{\epsilon}(0)$ and $\sigma_{\epsilon}(0)$. We introduce the notations
\be
X'_{\epsilon} = \fr{dX_{\epsilon}}{d\epsilon},
\la{xprime}
\ee 
with Eulerian form
\be
\eta_{\epsilon} = X'_{\epsilon}\circ X^{-1}_{\epsilon},
\la{eta}
\ee
also
\be
\tau'_{\epsilon} = \fr{d\tau_{\epsilon}}{d\epsilon},
\la{tauprime}
\ee
with Eulerian form
\be
\delta_{\epsilon} = \tau'_{\epsilon}\circ X^{-1}_{\epsilon},
\la{delta}
\ee
and
\be
u'_{\epsilon, 0} = \fr{du_{\epsilon}(0)}{d\epsilon}.
\la{uprimezero}
\ee
Differentiating ${\mathcal U}$ in (\ref{lagsys}) with respect to epsilon results in
\be
{\mathcal U}' = ({\mathbb L}(\na u_{\epsilon}(0))\circ X_{\epsilon})X'_{\epsilon}+{\mathbb L}(u'_{\epsilon, 0})\circ X_{\epsilon} + \left(\left(\na{\mathbb U}\right)(\sigma_{\epsilon})\circ X_{\epsilon}\right)X'_{\epsilon}-{\mathbb{U}}(\eta_{\epsilon}\cdot\na\sigma_{\epsilon})\circ X_{\epsilon} + {\mathbb U}(\delta_{\epsilon})\circ X_{\epsilon}.
\la{xprimeeq}
\ee
This follows from the fact that, for fixed $\tau$, 
\[
\fr{d}{d\epsilon}(\tau\circ X_{\epsilon}^{-1}) = -\na(\tau\circ X_{\epsilon}^{-1})(X'_{\epsilon}\circ X_{\epsilon}^{-1}),
\]
which in turn follows from
\[
\fr{dX_{\epsilon}^{-1}}{d\epsilon} = -(\na X_{\epsilon}^{-1})(X'_{\epsilon}\circ X_{\epsilon}^{-1}).
\]
Composing with $X_{\epsilon}^{-1}$ from the right, and dropping epsilon for ease of notation, we deduce from (\ref{xprimeeq})
\be
{\mathcal U}'\circ X^{-1} = {\mathbb L}(\na u_0)\eta + {\mathbb L}(u'_0) +
\left[\eta\cdot \na, \; {\mathbb U}\right](\sigma) +
{\mathbb U}(\delta)
\la{xprimeeul}
\ee
Here 
\be
[\eta\cdot\na,\; {\mathbb U}](\sigma) = \eta(t)\cdot\na{\mathbb U}(\sigma) - {\mathbb U}(\eta(s))\cdot\na \sigma)
\la{comm}
\ee
is the space-time commutator.  Note that
\be
X'(0) = \eta(0) = 0,\quad \delta(0) = \tau'_0 =\sigma'_0
\la{xidtauid}
\ee
and therefore
\be
{\mathcal U}'(0) = u'_0.
\la{uprimeid}
\ee
Differentiating ${\mathcal T}$ in (\ref{lagsys}) with respect to $\epsilon$ we obtain
\be
{\mathcal T}' = D_1F(g,\tau)g' + D_2F(g,\tau)\tau'
\la{tauprimeeq}
\ee
where
\be
g' = ({\mathbb L}(\na\na u_0)\circ X )X' + {\mathbb L}(\na u'_0)\circ X  + [\na {\mathbb G}(\sigma)\circ X]X' - {\mathbb G}(\eta\cdot \na\sigma)\circ X  + {\mathbb G}(\delta)\circ X
\la{gprime}
\ee
is the epsilon derivative of
\be 
g = {\mathbb L}(\na u_0)\circ X + {\mathbb{G}}( \tau\circ X^{-1})\circ X.
\la{gdef}
\ee
Composing with $X^{-1}$ we obtain
\be
g'\circ X^{-1} = {\mathbb L}(\na\na u_0)\eta + {\mathbb L}(\na u'_0) + [\eta \cdot\na,\; {\mathbb G}](\sigma) + {\mathbb G}(\delta)
\la{gprimeeul}
\ee
where
\be
[\eta\cdot\na, {\mathbb G}](\sigma) =\eta(t)\cdot\na {\mathbb G}(\sigma)-
{\mathbb G}(\eta(s)\cdot\na \sigma).
\la{gcomm}
\ee
Summarizing we have
\be
\left\{
\ba
\mathcal U'\circ X^{-1} = {\mathbb L}(\na u_0)\eta + {\mathbb L}(u'_0) + [\eta\cdot\na, {\mathbb U}](\sigma) + {\mathbb U}(\delta)\\
\mathcal T' = D_1F(g,\tau)g' + D_2F(g,\tau)\tau'\\ 
g'\circ X^{-1}= {\mathbb L}(\na\na u_0)\eta + {\mathbb L}(\na u'_0) + [\eta\cdot\na, {\mathbb G}](\sigma) + {\mathbb G}(\delta)
\ea
\right.
\la{linlageq}
\ee
where
\be
\eta = X'\circ X^{-1},\quad \delta = \tau'\circ X^{-1}.
\la{etadelta}
\ee
Differentiating $\mathcal U$ with respect to the Lagrangian independent variable (label) $a$ we have
\be
(\na_a {\mathcal U})(a,t) = \left({\mathbb L}(\na u_0)\circ X + {\mathbb{G}}( \tau\circ X^{-1})\circ X\right)(\na X) = g(a,t)(\na X)(a,t)
\la{namathcalu}
\ee
and using the fact that $\fr{d}{d\epsilon}$ and label derivatives commute we have
\be
\na_a {\mathcal U}'(a,t) = g'(a,t)(\na_a X(a,t)) + g(a,t)(\na_a X'(a,t))
\la{naxiprimeeq}
\ee
with $g$ given in(\ref{gdef}) above and $g'$ given by (\ref{gprime}).

\section{Bounds on operators and commutators}
We consider function spaces
\[
C^{\alpha, p} = C^{\alpha}(\Rr^d)\cap L^p(\Rr^d)
\]
with norm
\[
\|f\|_{\alpha, p} = \|f\|_{C^{\alpha}(\Rr^d)} + \|f\|_{L^p(\Rr^d)}
\]
for $\alpha\in (0,1)$, $p\in (1,\infty)$,  $C^{1+\alpha}(\Rr^d)$ with norm
\[
\|f\|_{C^{1+\alpha}(\Rr^d)} = \|f\|_{L^{\infty}(\Rr^d)} + \|\na f\|_{C^{\alpha}(\Rr^d)}
\]
and
\[ 
C^{1+\alpha, p} = C^{1+\alpha}(\Rr^d)\cap W^{1,p}(\Rr^d)
\]
with norm
\[
\|f\|_{1+\alpha, p} = \|f\|_{C^{1+\alpha}(\Rr^d)} + \|f\|_{W^{1,p}(\Rr^d)}
\]
We need also spaces of paths,
$L^{\infty}(0,T; C^{\alpha, p})$ with the usual norm, 
\[
\|f\|_{L^{\infty}(0,T; C^{\alpha,p})} = \sup_{t\in [0,T]}\|f(t)\|_{\alpha,p},
\]
spaces $Lip(0,T; C^{\alpha,p})$ with norm
\[
\|f\|_{Lip(0,T; C^{\alpha,p})} = \sup_{t\neq s, t,s\in [0,T]}\fr{\|f(t)-f(s)\|_{\alpha,p}}{|t-s|} + \|f\|_{L^{\infty}(0,T; C^{\alpha, p})},
\]
$C^{\beta}(0,T; C^{\alpha, p})$ with norm
\[
\|f\|_{C^{\beta}(0,T; C^{\alpha,p})} = \sup_{t\neq s, t,s\in [0,T]}\fr{\|f(t)-f(s)\|_{\alpha,p}}{|t-s|^{\beta}} + \|f\|_{L^{\infty}(0,T; C^{\alpha, p})}
\]
and spaces 
$C^{\beta}(0,T; C^{1+\alpha,p})$ with norm
\[
\|f\|_{C^{\beta}(0,T; C^{1+\alpha,p})} = \sup_{t\neq s, t,s\in [0,T]}\fr{\|f(t)-f(s)\|_{C^{1+\alpha, p}}}{|t-s|^{\beta}} + \|f\|_{L^{\infty}(0,T; C^{1+\alpha, p})}.
\]

We start with bounds on ${\mathbb U}$ and ${\mathbb G}$.
\beg{thm}\la{uboun} Let $0<\alpha<1$, $1<p<\infty$  and let $T>0$. There exists a constant such that
\be
\|\mathbb U(\sigma)\|_{L^{\infty}(0,T; C^{\alpha,p})}\le C\sqrt{T}\|\sigma\|_{L^{\infty}(0,T; C^{\alpha, p})}
\la{usbound}
\ee
and
\be
\|\mathbb{L}(u_0)\|_{L^{\infty}(0,T; C^{\alpha, p})} \le C\|u_0\|_{\alpha,p}
\la{lbound}
\ee
hold.
\end{thm} 

\beg{thm}\la{gbound} Let $0<\alpha<1$, $0<\beta\le 1$, $1<p<\infty$ and let $T>0$. The linear operator
\[
\sigma \mapsto {\mathbb G}(\sigma)
\]
\[
{\mathbb G}: C^{\beta}(0, T; C^{\alpha, p})\to
L^{\infty}(0,T; C^{\alpha,p})
\]
is continuous.
The linear operator 
\[
u_0\mapsto {\mathbb L}(\na u_0)
\]
maps continuously $C^{1+\alpha}(\Rr^d)\cap W^{1,p}(\Rr^d)$ to $L^{\infty}(0,T; C^{\alpha, p})$. 
\end{thm}
\beg{rem}
${\mathbb G}$ is actually continuous with values in the Banach space 
$C^{\beta\log}(0,T; C^{\alpha}(\Rr^d)\cap L^p(\Rr^d))$ with norm
\[
\ba
\|\gamma\|_{C^{\beta\log}(C^{\alpha,p})} = \sup_{t\in[0,T]}\|\gamma(\cdot, t)\|_{{\alpha,p}} + \\
\sup_{t,s\in [0,T], t\neq s}\fr{\|\gamma(\cdot, t)-\gamma(\cdot,s)\|_{{\alpha, p}}}{|t-s|^{\beta}\log {\left(\fr{2T}{|t-s|}\right)} + \log\left (1+\fr{|t-s|}{t}\right )}.
\ea
\]
Note that there is an unavoidable singularity at $t=0$ because if $\sigma\in C^{\alpha, p}$ is time-independent, then 
\[
{\mathbb G}(\sigma) = (\mathbb I-e^{t\Delta})R\mathbb H R\cdot\sigma
\] 
is bounded, but not time-H\"{o}lder continuous in $C^{\alpha, p}$ at $t=0$.
\end{rem}
\beg{thm}\la{commthm} Let $0<\alpha<1$, $\fr{1}{2}<\beta\le 1$, $T>0$. The bilinear operator
\[
(\eta,\sigma)\mapsto [\eta\cdot\na, {\mathbb G}](\sigma)
\]
is continuous from
\[
(\eta, \sigma)\in C^{\beta}(0,T; C^{1+\alpha}(\Rr^d) \times C^{\beta}(0,T; C^{\alpha, p})
\]
to 
\[
L^{\infty}(0,T; C^{\alpha, p}),
\]
that is,
\be
\|[\eta\cdot\na, {\mathbb G}]\sigma\|_{L^{\infty}(0,T; C^{\alpha,p})}\le
C\|\eta\|_{C^{\beta}(0,T; C^{1+\alpha}(\Rr^d))}\|\sigma\|_{C^{\beta}(0,T; C^{\alpha,p})}.
\la{bil}
\ee
\end{thm}
\beg{thm}\la{ucomm} Let $0<\alpha<1$, $0<\beta\le 1$, $T>0$. The bilinear operator
\[
(\eta, \sigma)\mapsto [\eta\cdot\na, {\mathbb U}](\sigma)
\]
obeys
\be
\ba
\|[\eta\cdot\na, {\mathbb U}](\sigma)\|_{L^{\infty}(0,T; C^{\alpha, p})}\\\le C[T^{1-\beta}\|\eta\|_{C^{\beta}(0,T; C^{\alpha}(\Rr^d))} + T^{\fr{1}{2}}\|\eta\|_{L^{\infty}(0,T; C^{1+\alpha}(\Rr^d))}]\|\sigma\|_{L^{\infty}(0,T; C^{\alpha, p})} 
\ea
\la{ucommbound}
\ee
\end{thm}

\noindent{\bf{Proof of Theorem  {\ref{gbound}}}}. \; We note first that 
\[
u_0\mapsto e^{t\Delta}\na u_0 = {\mathbb L}(\na u_0)
\]
maps continuously $C^{1+\alpha}(\Rr^d)\cap W^{1,p}(\Rr^d)$ to $L^{\infty}(0,T; C^{\alpha, p})$. We consider therefore  ${\mathbb G}(\sigma)$. 
We write
\[
{\mathbb G}(\sigma)(t) = \int_0^te^{(t-s)\Delta}\na{\mathbb H}\na \sigma(s)ds
=G_1 + G_2
\]
with
\[
G_1(x,t) = \int_0^{t-l}e^{(t-s)\Delta}\na{\mathbb H}\na \sigma(s)ds
\]
and 
\[
G_2(x,t) = \int_{t-l}^te^{(t-s)\Delta}\na{\mathbb H}\na \sigma(s)ds
\]
and $0<l<t$ arbitrary, to be chosen later. The solution of the heat equation is given by convolution with the Gaussian $g_t$. We use the fact that
$\na g_{t-s}$ is in the Hardy class $H^1(\Rr^d)$ and therefore
\[
\|\na\mathbb H\na g_{t-s}\|_{L^1(\Rr^d)}\le \fr{C}{t-s}
\]
to deduce
\[
\ba
\|G_1(\cdot, t)\|_{\alpha, p} \le C\int_0^{t-l}\|\na{\mathbb H}\na g_{t-s}\|_{L^1(\Rr^d)}\|\sigma (s)\|_{\alpha,p}ds\\
\le C\|\sigma\|_{L^{\infty}(0,T; C^{\alpha,p})}\log\left(\fr{t-l}{t}\right)
\ea
\]
We split $G_2$ in two pieces
\[
G_2(x,t) = G_3 + G_4
\]
with 
\[
G_3(x,t) = \int_{t-l}^te^{(t-s)\Delta}\na{\mathbb H}\na (\sigma(s)-\sigma(t))ds
\]
and
\[
G_4(x,t) = \int_{t-l}^te^{(t-s)\Delta}\na{\mathbb H}\na \sigma(t)ds.
\]
Now
\[
\|G_3(\cdot, t)\|_{\alpha, p}\le C\|\sigma\|_{C^{\beta}(C^{\alpha,p})}\int_{t-l}^t(t-s)^{-1+\beta}ds = C\|\sigma\|_{C^{\beta}(C^{\alpha,p})} l^{1-\beta}
\]
and, because
\[
G_4(\cdot, t) = \left({\mathbb{I}} - R{\mathbb H}Re^{l\Delta}\right)e^{(t-l)\Delta}
\]
we obtain
\[
\|G_4(\cdot, t)\|_{C^{\alpha,p}}\le C\|\sigma\|_{L^{\infty}(0,T; C^{\alpha,p})}.
\]
Choosing $l=\fr{t}{2}$ we obtain
\[
\|{\mathbb G}\sigma\|_{L^{\infty}(0,T; C^{\alpha,\beta})}\le C(1+ T^{1-\beta})
\|\sigma\|_{C^{\beta}(0,T; C^{\alpha, p})}
\]
and this ends the proof of the theorem. 

Considering the remark following the theorem, let $0\le t<t+h\le T$, and let us write
\[
\ba
{\mathbb G}(\sigma)(t+h)-{\mathbb G}(\sigma)(t) =\\ 
\int_0^te^{(t-s)\Delta}\na\mathbb H\na(\sigma(s+h)-\sigma(s))ds +
\int_0^he^{(t+h-s)\Delta}\na{\mathbb H}\na\sigma(s)ds=\\
I_1(x,t) + I_2(x,t).
\ea
\]
We write 
\[
\ba
I_1(x,t) = I_{11}(x,t) + I_{12}(x,t) = \int_0^te^{(t-s)\Delta}\na\mathbb H\na(\sigma(t+h)-\sigma(t))ds +\\
\int_0^te^{(t-s)\Delta}\na\mathbb H\na(\sigma(s+h)-\sigma(t+h) +\sigma(t)-\sigma(s))ds
\ea
\]
Because
\[
\int_0^te^{(t-s)\Delta}\na\mathbb H\na ds = (\mathbb I-e^{t\Delta})R\mathbb H R
\]
we obtain
\[
\|I_{11}(\cdot,t)\|_{C^{\alpha}} \le Ch^{\beta}\|\sigma\|_{C^{\beta}(0,T; C^{\alpha}(\Rr^d))}.
\]
A logarithm is lost in the estimate of $I_{12}$. Using the properties of $g_{t-s}$, we have 
\[
\ba
\|I_{12}(\cdot,t)\|_{C^{\alpha}(\Rr^d)}\le \\
C\|\sigma\|_{C^{\beta}(0,T; C^{\alpha}(\Rr^d))}\int_0^{t}\fr{1}{t-s}\min\{|t-s|^{\beta}; h^{\beta}\}ds\le \\
C\|\sigma\|_{C^{\beta}(0,T; C^{\alpha}(\Rr^d))}h^{\beta}(1+ \log_{+}{\fr{t}{h}})\le 
C\|\sigma\|_{C^{\beta}(0,T; C^{\alpha}(\Rr^d))}h^{\beta}(1+ \log{\fr{2T}{h}})
\ea
\]
Thus, for $I_1$ we obtain uniform H\"{o}lder continuity in time on $[0,T]$ with logarithmic loss, and in particular with any exponent less than $\beta$.
For $I_2$ we have
\[
\ba
\|I_2\|_{C^{\alpha}(\Rr^d)}\le C\int_0^h\|\na{\mathbb H}\na g_{t+h-s}\|_{L^1(\Rr^d)}\|\sigma(\cdot,s)\|_{C^{\alpha}(\Rr^d)}ds\le\\
C\|\sigma\|_{L^{\infty}(0,T; C^{\alpha}(\Rr^d))}\int_0^h\fr{1}{t+h-s}ds =
C\|\sigma\|_{L^{\infty}(0,T; C^{\alpha}(\Rr^d))}\log\left(1+\fr{h}{t}\right).
\ea
\]
Thus $I_2$ is Lipschitz continuous in time away from $t=0$, but not H\"{o}lder continuous in time at $t=0$. Because $\log(1 +x)\le Cx^{\beta}$, if we measure H\"{o}lder continuity in time of order $\beta$ we have a
singular coefficient of order $t^{-\beta}$ near $t=0$.

\beg{lemma}\la{steadycomm} Let $0<\alpha<1$, $1<p<\infty$. Let $\eta\in C^{1+\alpha}(\Rr^d)$ and let
\[
({\mathbb K}\sigma)(x) = P.V.\int_{\Rr^d}k(x-y)\sigma(y)dy
\]
be a classical Calderon-Zygmund operator with kernel $k$ which is smooth away from the origin, homogeneous of degree $-d$ and with mean zero on spheres about the origin. Then the commutator  $[\eta\cdot\na, {\mathbb K}]$ can be defined as a bounded linear operator in $C^{\alpha, p}$ 
and
\be
\|[\eta\cdot\na, {\mathbb K}]\sigma\|_{\alpha, p}\le C\|\eta\|_{C^{1+\alpha}(\Rr^d)}\|\sigma\|_{\alpha, p}
\la{bilsteady}
\ee
\end{lemma}
\beg{rem}
The conclusion of the lemma holds also for operators ${\mathbb H}$ which are products of classical CZ operators. This follows from telescoping applications of the lemma.
\end{rem}
\noindent{\bf{Proof of Lemma (\ref{steadycomm})}}. Let us note that both terms in the commutator, $\eta\cdot\na {\mathbb K}\sigma$ and $\mathbb K\eta\na\sigma$ are well defined and H\"{o}lder continuous if $\sigma $ is smooth. We compute first
\[
[\eta\cdot\na, {\mathbb K}]\sigma(x) = \int_{\Rr^d}k(x-y)(\eta(x)-\eta(y))\cdot\na_y\sigma(y)dy
\]
Now we introduce a smooth cutoff $\chi(|x-y|)$ identically equal to 1 for $|x-y|\le 1$ and compactly supported. The conclusion of the lemma holds for
\[
C_{out}(x) = \int_{\Rr^d}\left(1-\chi(|x-y|)\right)k(x-y)(\eta(x)-\eta(y))\cdot\na_y\sigma(y)dy
\]
by integration by parts and inspection, using the $L^p$ bound for $\sigma$. 

We concentrate our attention on
\[
C_{in}(x) = \int_{\Rr^d}\chi(|x-y|)k(x-y)(\eta(x)-\eta(y))\cdot\na_y\sigma(y)dy
\]
We first write
\[
C_{in}(x) = \int_{\Rr^d}\chi(|x-y|)k(x-y)(\eta(x)-\eta(y))\cdot\na_y(\sigma(y)-\sigma(x))dy
\]
and then we integrate by parts:
\[
C_{in}(x) = C(x) + C_1(x) + C_2(x)
\]
with 
\[
C(x) = \int_{\Rr^d}\chi(|x-y|)\na_xk(x-y)(\eta(x)-\eta(y))(\sigma(y)-\sigma(x))dy,
\]
\[
C_1(x) = \int_{\Rr^d}\na_x\chi(|x-y|))k(x-y)(\eta(x)-\eta(y))(\sigma(y)-\sigma(x))dy
\]
and
\[
C_2(x) = \int_{\Rr^d}\chi(|x-y|)k(x-y)(\na_y\eta(y))(\sigma(y)-\sigma(x))dy.
\]
It is easy to see that $C_1$ and $C_2$ are H\"{o}lder continuous and satisfy the bound (\ref{bilsteady}). We investigate now $C(x)$ and write
\[
\eta(y)-\eta(x) = (y-x)\cdot\int_0^1 \na\eta(x+\lambda(y-x))d\lambda
\]
We write also
\[
K(z) = z\na_zk (z)
\]
and note that it is homogeneous of order $-d$ and smooth away from the origin. The averages on spheres might not vanish.
So, with these preparations $C(x)$ is
\[
C(x) = \int_0^1d\lambda \int_{\Rr^d}\chi(|x-y|)K(x-y)(\na\eta(x+\lambda(y-x)))(\sigma(y)-\sigma(x))dy
\]
We write now
\[
C(x) = A(x) + B(x)
\]
where
\[
A(x) = \int_0^1d\lambda \int_{\Rr^d}\chi(|x-y|)K(x-y)(\na\eta(x+\lambda(y-x)) -\na\eta(x))(\sigma(y)-\sigma(x))dy
\]
and
\[
B(x) =\na\eta(x)\int_{\Rr^d}\chi(|x-y|)K(x-y)(\sigma(y)-\sigma(x))dy.
\]
Now $B\in C^{\alpha}(\Rr^d)$ and obeys (\ref{bilsteady}). It is obviously enough to check that
\[
I(x) = \int_{\Rr^d}\chi(|x-y|)K(x-y)(\sigma(y)-\sigma(x))dy
\]
is in $C^{\alpha}(\Rr^d)$ and its norm is bounded by that of $\sigma$. To check this we take the difference
\[
I(x+h)-I(x) = I_1+I_2 +I_3
\]
where $I_1$ is 
\[
\ba
I_1 = \int_{|x-y|\le 4|h|,\; |x+h-y|\ge 4|h|}\chi(|x+h-y|)K(x+h-y)(\sigma(y)-\sigma(x+h))dy\\
-\int_{|x-y|\le 4|h|, \; |x+h-y|\ge 4|h|}\chi(|x-y|)K(x-y)(\sigma(y)-\sigma(x))dy
\ea
\]
$I_{2}$ is 
\[
\ba
I_2 = \int_{|x-y|\ge 4|h|,\; |x+h-y|\le 4|h|}\chi(|x+h-y|)K(x+h-y)(\sigma(y)-\sigma(x+h))dy\\
-\int_{|x-y|\ge 4|h|,\; |x+h-y|\le 4|h|}\chi(|x-y|)K(x-y)(\sigma(y)-\sigma(x))dy
\ea
\]
and
\[
\ba
I_3 = \int_{|x-y|\ge 4|h|,\; |x+h-y|\ge 4|h|}\chi(|x+h-y|)K(x+h-y)(\sigma(y)-\sigma(x+h))dy\\
-\int_{|x-y|\ge 4|h|,\; |x+h-y|\ge 4|h|}\chi(|x-y|)K(x-y)(\sigma(y)-\sigma(x))dy
\ea
\]
For $I_1$ and $I_2$ we note that both $|x-y|\le 5|h|$ and $|x+h-y|\le 5|h|$, and we use the straightforward inequality
\[
\int_0^{5h}r^{-1}r^{\alpha}dr \le Ch^{\alpha}.
\]
The integral $I_3$ is split into two pieces.
\[
I_3 = I_4 + I_5
\]
with
\[
I_4 =\fr{1}{2}\int_{|x-y|\ge 4|h|,\; |x+h-y|\ge 4|h|}[\chi(|x+h-y|)K(x+h-y)-\chi(|x-y|)K(x-y)](2\sigma(y)-\sigma(x+h)+\sigma(x))dy
\]
and 
\[
I_5 = \fr{\sigma(x+h)-\sigma(x)}{2}\int_{|x-y|\ge 4|h|,\; |x+h-y|\ge 4|h|}[\chi(|x+h-y|)K(x+h-y)+\chi(|x-y|)K(x-y)]dy.
\]
For $I_4$ we use the smoothness of the kernel, the intermediate value theorem,
and the H\"{o}lder bounds to obtain
\[
|I_4|\le C\|\sigma\|_{C^{\alpha}(\Rr^d)}\int_{3|h|}^{\infty}|h|r^{-2}(r^{\alpha}+ h^{\alpha})dr\le C\|\sigma\|_{C^{\alpha}(\Rr^d)}|h|^{\alpha}.
\]
For $I_5$ we recall that $K(z) = z\na k(z)$. We claim that integrals of
\[
\int_{|z|\ge 4|h|, |z\pm h|\ge 4|h|}\chi(|z|)K(z)dz
\]
are bounded uniformly, independently of $h$. Indeed, integrating by parts
\[
\ba
\int_{|z|\ge 5|h|}\chi(|z|)z_j\pa_i k(z)dz =\\ -\int_{|z|\ge 5|h|}\delta_{ij}\chi(|z|)k(z)dz -\int_{|z|\ge 5|h|}\pa_i(\chi(|z|))z_jk(z)dz + \int_{|z|=5|h|}z_j\chi(5|h|)\fr{z_i}{5|h|}k(z)dS(z)\\
= 0 +{\mbox{ bounded}}.
\ea
\]
Here we used that $k$ has mean zero on spheres. On the other hand, on the annular regions we use simply the homogeneity of $K$ and
\[
\int_{4h}^{5h}\fr{1}{r}dr \le C.
\]
The integral $A(x)$ is treated in a similar fashion. We write
\[
A(x) = \int_{\Rr^d}\chi(|x-y|)K(x-y)\phi(x,x-y)(\sigma(y)-\sigma(x))dy
\]
where
\[
\phi(x,x-y) = \int_0^1(\na\eta(x+\lambda(y-x))-\na\eta(x))d\lambda
\]
We consider 
\[
A(x+h)-A(x) = A_1 + A_2 +A_3
\]
where $A_1$ and $A_2$, like  $I_1$ and $I_2$ above, are differences of integrals on $|x-y|\le 4|h|$ and $|x+h-y|\ge 4|h|$, and, respectively, on $|x-y|\ge 4|h|$ and $|x+h-y|\le 4|h|$, while $A_3$ is the difference of integrals corresponding to both
$|x-y|\ge 4|h|$ and $|x+h-y|\ge 4|h|$. As before, using the triangle inequality, the regions of integration for $A_1$ and $A_2$ are regions where both $|x-y|\le 5|h|$ and $|x+h-y|\le 5|h|$ and  therefore, the integrals are small separately, without need to take the difference. Using the fact that
\[
|\phi(x,x-y))|\le \|\eta\|_{C^{1+\alpha}(\Rr^d)}|x-y|^{\alpha}
\]
we obtain that
\[
|A_1| + |A_2| \le C|h|^{2\alpha}\|\sigma\|_{C^{\alpha}(\Rr^d)}\|\eta\|_{C^{1+\alpha}(\Rr^d)}
\]
We treat $A_3$ as we treated $I_3$: we split $A_3 = A_4 +A_5$, where 
\[
\ba
A_4 =\fr{1}{2}\int_{|x-y|\ge 4|h|,\; |x+h-y|\ge 4|h|}[\chi(|x+h-y|)K(x+h-y)-\chi(|x-y|)K(x-y)]\\\times[\phi(x+h,x+h-y)(\sigma(y)-\sigma(x+h)) + \phi(x,x-y)(\sigma(y)-\sigma(x))]dy
\ea
\]
and, using the smoothness of the kernel and the bounds on $\phi$ and $\sigma$, this leads to an integral inequality 
\[
\int_{3h}^1hr^{-2}(r^\alpha +h^{\alpha})^2\le Ch^{2\alpha}
\]
so
\[
|A_4|\le  C|h|^{2\alpha}\|\sigma\|_{C^{\alpha}(\Rr^d)}\|\eta\|_{C^{1+\alpha}(\Rr^d)}
\]
Finally we treat
\[
\ba
A_5 = \fr{1}{2}\int_{|x-y|\ge 4|h|,\; |x+h-y|\ge 4|h|}[\chi(|x+h-y|)K(x+h-y)+\chi(|x-y|)K(x-y)]\\
\times[\phi(x+h,x+h-y)(\sigma(y)-\sigma(x+h))- \phi(x,x-y)(\sigma(y)-\sigma(x))]dy.
\ea
\]
We note by polarization that
\[
\ba
\left|\phi(x+h,x+h-y)(\sigma(y)-\sigma(x+h))- \phi(x,x-y)(\sigma(y)-\sigma(x))\
\right| \\
\le C\|\sigma\|_{C^{\alpha}(\Rr^d)}\|\eta\|_{C^{1+\alpha}(\Rr^d)}(|x-y| +|h|)^{\alpha}|h|^{\alpha}
\ea
\]
and therefore $A_5$ is bounded directly using
\[
\int_{4h}^1r^{-1}(h^{\alpha}(r^{\alpha}+h^{\alpha})dr\le Ch^{\alpha}
\]
\[
|A_5|\le C\|\sigma\|_{C^{\alpha}(\Rr^d)}\|\eta\|_{C^{1+\alpha}(\Rr^d)}|h|^{\alpha}.
\]
This concludes the proof of the fact that
\[
\|C_{in}\|_{C^{\alpha}(\Rr^d)} \le C\|\sigma\|_{\alpha, p}\|\eta\|_{C^{1+\alpha}(\Rr^d)}
\]
The proof of the $L^p$ bound in Lemma (\ref{steadycomm}) is done using the observation that
\[
C(\sigma) = [\eta\cdot\na, {\mathbb K}]\sigma = PV\int_{\Rr^d} K(x,y)\sigma(y)dy  - {\mathbb K}((\na\cdot\eta)\sigma) 
\]
where 
\[
K(x,y) = (\eta_j(x)-\eta_j(y)))\pa_j k(x-y).
\]
Now the operator ${\mathbb K}$ is bounded in $L^p$ spaces and the operator $T$ given by
\[
\sigma \mapsto PV\int_{\Rr^d} K(x,y)\sigma(y)dy = (T\sigma)(x)
\]
is a Calderon-Zygmund operator, that is, the kernel $K$ is smooth away from the diagonal, obeys
\[
|K(x,y)|\le C\fr{1}{|x-y|^d}
\]
and
\[
|K(x+h,y)-K(x,y)| + |K(x,y+h)-K(x,y)| \le C\fr{|h|}{|x-y|^{d+1}}
\]
for $2|h| \le |x-y|$, and  $T$ is bounded in $L^2(\Rr^d)$. The boundedness in $L^2$ is verified quickly below. It follows that $T$ is bounded in $L^p(\Rr^d)$, $1<p<\infty$ (see for instance  cite: stein). 
For the bound in $L^2$ we need to verify that
\[
\left|\int_{\Rr^d}(Tf)(x)g(x)dx\right| \le C\|f\|_{L^2(\Rr^d)}\|g\|_{L^2(\Rr^d)}.
\]
We write
\[
\ba
\int_{\Rr^d}(Tf)(x)g(x)dx =\int_{\Rr^d}dx PV\int_{|z|\le 1} K(x,x-z)f(x-z)g(x) dz \\ + \int_{\Rr^d}dx \int_{|z|\ge 1}K(x,x-z)f(x-z)g(x)dz
= T_1 + T_2
\ea
\]
Clearly
\[
|T_2| \le C\|\eta\|_{L^{\infty}(\Rr^d)}\|\|f\|_{L^2(\Rr^d)}\|g\|_{L^2(\Rr^d)}
\]
because
\[
|K(x,y)|\le C\|\eta\|_{L^{\infty}(\Rr^d)} \fr{1}{|x-y|^{d+1}}
\]
in view of the homogeneity of $k$. For $T_1$ we use the fact that we have
\[
\left| PV\int_{|z|\le 1}(\eta(x)-\eta(x-z))\na k(z)dz\right | \le C\|\eta\|_{C^{1+\alpha}}
\]
uniformly in $x\in \Rr^d$. Indeed, this is easily verified in a manner similar to the proof in the $C^{\alpha}$ case, using the fact that
\[
\eta(x) - \eta(x-z) - z\na\eta(x) = O(|z|^{1+\alpha}),
\]
integration by parts, and the vanishing of spherical averages of $k$:
\[
\ba
\left |PV\int_{|z|\le 1}(\eta(x)-\eta(x-z))\na k(z)dz\right| \le \int_{|z|\le 1}|(\eta(x)-\eta(x-z)-z\na\eta(x))||\na k(z)|dz\\
+|\na\eta(x)|\left|PV\int_{|z|\le 1} z\na k(z)dz\right|\le C\|\eta\|_{C^{1+\alpha}(\Rr^d)}.
\ea
\]

\noindent{\bf{Proof of Theorem \ref{commthm}}}.

The commutator is 
\[
C(x,t)= [\eta\cdot\na, {\mathbb G}](\sigma) =
\eta(\cdot,t)\cdot\na\int_0^te^{(t-s)\Delta}\na\mathbb H\na\sigma(s)ds -
\int_0^te^{(t-s)\Delta}\na\mathbb H\na(\eta(s)\cdot\na\sigma(s))ds.
\]
We need to show that $\|C\|_{L^{\infty}(0, T; C^{\alpha}(\Rr^d))}$ is finite, under our assumptions on $\eta$ and $\sigma$.
Recall our notation that $e^{(t-s)\Delta}$ is given by convolution with $g_{t-s}$. We start by the observation that we can replace $\eta(s)$ by $\eta(t)$ in the second term of the commutator. Indeed 
\[
\int_0^te^{(t-s)\Delta}\na\mathbb H\na((\eta(s)-\eta(t))\cdot\na\sigma(s))ds = I_1(t) + I_2(t)
\]
with 
\[
I_1(t) =-\int_0^te^{(t-s)\Delta}\na\mathbb H\na\{[\na\cdot(\eta(s)-\eta(t))]\sigma(s)\}ds
\]
and
\[
I_2(t) = -\int_0^t(\na\mathbb H\na\na g_{t-s})*\{(\eta(s)-\eta(t))\sigma(s)\}ds.
\]
Now
\[
I_1 = \mathbb G\{[\na\cdot(\eta(s)-\eta(t))]\sigma\}
\]
and it belongs to $L^{\infty}(0, T; C^{\alpha, p})$ by Theorem {\ref{gbound}} because $(\na\cdot\eta)\in C^{\beta}(0,T; C^{\alpha}(\Rr^d))$. 
Regarding $I_2$, we have
\[
\ba
\|I_2(t)\|_{\alpha, p}\le C\int_0^t\|\na\mathbb H\na\na g_{t-s}\|_{L^1(\Rr^d)}\|(\eta(s)-\eta(t))\sigma(s))\|_{\alpha, p}ds\\
\le C\|\eta\|_{C^{\beta}(0,T; C^{\alpha}(\Rr^d))}\|\sigma\|_{L^{\infty}(0, T; C^{\alpha, p})}\int_0^t (t-s)^{-\fr{3}{2}}(t-s)^{\beta}ds\\
\le Ct^{\beta-\fr{1}{2}}\|\eta\|_{C^{\beta}(0,T; C^{\alpha}(\Rr^d))}\|\sigma\|_{L^{\infty}(0, T; C^{\alpha, p})}
\ea
\]
So now we have to examine
\[
C_1 = \eta(\cdot,t)\cdot\na\int_0^te^{(t-s)\Delta}\na\mathbb H\na\sigma(s)ds -\int_0^te^{(t-s)\Delta}\na\mathbb H\na(\eta(t)\cdot\na\sigma(s))ds.
\]
Now we observe that we can replace $\sigma(s)$ by $\sigma(s)-\sigma(t)$ in both integrals. Indeed, replacing $\sigma(s)$ in $C_1$ by $\sigma(t)$ integrates in time to
\[
\eta(\cdot,t)\cdot\na\int_0^te^{(t-s)\Delta}\na\mathbb H\na\sigma(t)ds -\int_0^te^{(t-s)\Delta}\na\mathbb H\na(\eta(t)\cdot\na\sigma(t))ds=
[\eta(t)\cdot\na, (\mathbb I-e^{t\Delta})R\mathbb HR]\sigma(t)
\]
and the commutator
\[
[\eta(t)\cdot\na, (\mathbb I-e^{t\Delta})R\mathbb HR]\sigma(t)
\]
is bounded by Lemma {\ref{steadycomm}}. It remains to investigate 
\[
\ba
C_2 = \eta(\cdot,t)\cdot\na\int_0^te^{(t-s)\Delta}\na\mathbb H\na(\sigma(s)-\sigma(t))ds -\int_0^te^{(t-s)\Delta}\na\mathbb H\na(\eta(t)\cdot\na(\sigma(s)-\sigma(t)))ds\\
= \eta(\cdot,t)\cdot\na\int_0^t(\na\na g_{t-s})*\mathbb H(\sigma(s)-\sigma(t))ds -\int_0^t(\na\na g_{t-s})*\mathbb H(\eta(t)\cdot\na(\sigma(s)-\sigma(t)))ds
\ea
\]
We claim that we can move $\eta\cdot\na$ inside the first time integral, past the convolution with the derivative of the heat kernel. Indeed,
the difference is
\[
D(x,t)= \int_0^t\int_{\Rr^d}(\na\na\na g_{t-s})(x-y)(\eta(x,t)-\eta(y,t)){\mathbb H}(\sigma(s)-\sigma(t))(y)dyds
\]
We use now the fact that
\[
\||z|\na\na\na g_{t-s}(z)\|_{L^{1}(\Rr^d)}\le C(t-s)^{-1}
\]
to deduce, after changing variables to $z=x-y$ and writing $\eta(x,t)-\eta(x-z,t) = -z\int_0^1\na\eta(x-\lambda z, t)d\lambda$ that
\[
\|D(\cdot, t)\|_{{\alpha, p}}\le C\|\eta\|_{L^{\infty}(0,T; C^{1+\alpha}(\Rr^d))}\|\sigma\|_{C^{\beta}(0,T; C^{\alpha, p})}\int_0^t(t-s)^{-1+\beta}ds
\]
So, finally we arrived at
\[
C_3(t) = \int_0^t(\na\na g_{t-s})*[\eta(t)\cdot\na, {\mathbb H}](\sigma(s)-\sigma(t))ds.
\]
We bound this using Lemma {\ref{steadycomm}}:
\[
\|C_3(t)\|_{\alpha, p}\le C\|\eta\|_{L^{\infty}(0,T; C^{1+\alpha}(\Rr^d))}\|\sigma\|_{C^{\beta}(0,T; C^{\alpha, p})}\int_0^t(t-s)^{-1+\beta}ds
\]
This concludes the proof of Theorem {\ref{commthm}}.

{\bf{Proof of Theorem{\ref{uboun}.}}} We bound
\[
\|{\mathbb U}(\sigma)\|_{\alpha, p}\le
C\int_0^t\|\na g_{t-s}\|_{L^1(\Rr^d)}\|\sigma(s)\|_{\alpha, p}ds\le
C\|\sigma\|_{L^{\infty}(0,t; C^{\alpha, p})}\int_0^t\fr{1}{\sqrt{t-s}}ds
\]

{\bf{Proof of Theorem {\ref{ucomm}}}}. The computation concerns
\[
\eta(t)\cdot\na \int_0^t e^{(t-s)\Delta}{\mathbb H}\na\sigma(s)ds -
\int_0^t e^{(t-s)\Delta}{\mathbb H}\na(\eta(s)\cdot\na\sigma(s))ds 
\]
Replacing $\eta(s)$ by $\eta(t)$ in the second term in  the commutator, introduces
\[
E(t) = \int_0^te^{(t-s)\Delta}{\mathbb H}\na((\eta(s)-\eta(t))\cdot\na\sigma(s))ds = E_1(t) + E_2(t)
\]
where
\[
E_1(t) = -{\mathbb U}((\na\cdot(\eta(s)-\eta(t))\sigma(s))
\]
and 
\[
E_2(t) = -\int_0^t \na\na g_{t-s}*{\mathbb H}((\eta(s)-\eta(t))\sigma(s))ds.
\]
The first term is bounded by (\ref{usbound})
\[
\|E_1(t)\|_{\alpha, p}\le C\sqrt{t}\|\eta\|_{L^{\infty}(0,T; C^{1+\alpha}(\Rr^d))}\|\sigma\|_{L^{\infty}(0,T; C^{\alpha, p})}
\]
and $E_2(t)$ is bounded by
\[
\|E_2(t)\|_{\alpha, p} \le Ct^{1-\beta}\|\eta\|_{C^{\beta}(0,T; C^{\alpha}(\Rr^d))}\|\sigma\|_{L^{\infty}(0,T; C^{\alpha, p})}.
\]
We have to bound now
\[
V(t) = \eta(t)\cdot\na\int_0^te^{(t-s)\Delta}{\mathbb H}\na\sigma(s)ds -
\int_0^te^{(t-s)\Delta}{\mathbb H}\na ((\eta(t)\cdot\na\sigma(s)))ds
\]
We claim that we can put $\eta(t)\cdot\na$ inside the first time integral, past the convolution with the gradient of the heat kernel. Indeed, the difference 
\[
D(x,t) = \int_0^t\int_{\Rr^d}\na\na g_{t-s}(z)(\eta(x,t)-\eta(x-z,t)){\mathbb H}\sigma (x-z,s)dz
\]
can be bounded, after writing $\eta(x,t)-\eta(x-z,t) = -z\int_0^1\na\eta(x-\lambda z, t)d\lambda$, and using
\[
\||z|\na\na g_{t-s}\|_{L^1(\Rr^d)}\le C\fr{1}{\sqrt{t-s}}
\]
by
\[
\|D\|_{L^{\infty}(0,T; C^{\alpha,p})} \le C\sqrt T\|\eta\|_{L^{\infty}(0,T; C^{1+\alpha}(\Rr^d))}\|\sigma\|_{L^{\infty}(0,T; C^{\alpha, p})}
\]
We are left with
\[
E= \int_0^t\na e^{(t-s)\Delta}[\eta(t)\cdot\na, {\mathbb H}]\sigma(s)ds + \int_0^t\na e^{(t-s)\Delta} (\div\eta(t)){\mathbb H}(\sigma(s))ds 
\]
which we bound using Lemma \ref{steadycomm}
\[
\|E\|_{L^{\infty}(0,T; C^{\alpha,p})} \le C\sqrt{T}\|\eta\|_{L^{\infty}(0,T; C^{1+\alpha}(\Rr^d))}\|\sigma\|_{L^{\infty}(0,T; C^{\alpha, p})}
\]

Summing the bounds we conclude that (\ref{ucommbound}) holds.
\section{Bounds on solutions}

We start with a few kinematic observations. Let $u\in L^{\infty}(0,T; C^{1+\alpha}(\Rr^d))$ be a velocity. Then the Lagrangian maps $X(a,t)$ are
\[
X(a,t) = a + \chi(a,t)
\]
with $\chi\in Lip(0,T; C^{1+\alpha}(\Rr^d))$, $\chi(a,0)=0$. Moreover $X^{-1}(x,t) = A(x,t)$
obeys the transport equation
\[
\pa_t A + u\cdot\na A =0
\]
with $A(x,0) = x$ and $A(x,t) = x + \alpha(x,t)$ with $\alpha \in Lip(0,T; C^{1+\alpha}(\Rr^d))$, $\alpha(x,0)=0$. (Obviously, $\chi(a,t) + \alpha (X(a,t),t)= 0$). The inverse exists even if $u$ is not divergence-free. The gradients obey
\[
\|\na X^{-1}(t)\|_{L^{\infty}(\Rr^d)} \le \exp{\int_0^t \|\na u(s)\|_{L^{\infty}(\Rr^d)}ds}
\]
The same is true for the gradients $\na X$:
\[
\|\na X(t)\|_{L^{\infty}(\Rr^d)} \le \exp{\int_0^t \|\na u(s)\|_{L^{\infty}(\Rr^d)}ds}
\]
Because
\[
a-b = X^{-1}(X(a,t),t) - X^{-1}(X(b,t),t)
\]
it follows that
\[
|a-b|\le |X(a,t)-X(b,t)|\exp{\int_0^t \|\na u(s)\|_{L^{\infty}(\Rr^d)}ds},
\]
and because
\[
X(a,t)-X(b,t) = \int_0^1\fr{d}{d\mu}X((1-\mu)a + \mu b, t)d\mu
\]
it follows that
\[
|X(a,t)-X(b,t)|\le |a-b|\exp{\int_0^t \|\na u(s)\|_{L^{\infty}(\Rr^d)}ds}.
\]
Therefore we have the important and quite general chord-arc bound
\be
\lambda^{-1}\le \fr{|a-b|}{|X(a,t)-X(b,t)|} \le \lambda
\la{ca}
\ee
where
\be
\lambda (t) = \exp{\int_0^t \|\na u(s)\|_{L^{\infty}(\Rr^d)}ds}.
\la{lambda}
\ee
Because of the chord-arc bound it is possible and convenient to measure the size of the Lagrangian nonlinearities in H\"{o}lder spaces after composition with $X^{-1}$, i.e. in Eulerian variables. We consider the equation (\ref{tauode}) now. Let us note that from, general ODE theory, we have a priori bounds for short time. Using the same notation (\ref{g}) for $\na u\circ X$ we have
\be
\|\tau(t)\|_{L^{\infty}(\Rr^d)}\le K
\la{taufty}
\ee
for $t\le T$, where $K$ and $T$ depend on $\|\sigma_0\|_{L^{\infty}(\Rr^d)}$ and a bound on $\sup_{[0,T]}\|\na u\|_{L^{\infty}(\Rr^d)}$. T
and, consequently for $\sigma = \tau\circ X^{-1}$,
\be
\|\sigma(t)\|_{L^{\infty}(\Rr^d)}\le K
\la{sigfty}
\ee
Similarly
\be
\|\tau(t)\|_{L^{p}(\Rr^d)}\le K_p
\la{taup}
\ee
and consequently
\be
\|\sigma(t)\|_{L^{p}(\Rr^d)}\le K_p\lambda(t)^{\fr{d}{p}}.
\la{sigp}
\ee
The exponent of $\lambda$ in the last inequality is $0$ if we assume incompressibility. 
Also, by taking finite differences in Lagrangian variables 
\[
\delta_h\tau (a,t) = \tau (a+h,t)-\tau(a,t)
\]
we obtain 
\[
\fr{d}{dt}\delta_h\tau = F(g(a+h), \tau(a+h)) -F(g(a),\tau(a))
\]
and deduce, via
\[
\fr{d}{dt}|\delta_h\tau| \le \|D_2F(g,\tau)\|_{L^{\infty}(\Rr^d)}|\delta_h\tau| +
\|D_1F(g,\tau)\|_{L^{\infty}(\Rr^d)}|\delta_h g|
\]
that
\be
\|\tau(t)\|_{C^{\alpha}(\Rr^d)}\le C_{\alpha}
\la{taualpha}
\ee
with $C_{\alpha}$ depending on $\|\sigma_0\|_{C^{\alpha}(\Rr^d)}$ and 
$\sup_{[0,T]}\|\na u\|_{C^{\alpha}(\Rr^d)}$.
Passing to the Eulerian seminorm costs $\lambda^{\alpha}$:
\be
\|\sigma(t)\|_{C^{\alpha}(\Rr^d)}\le C_{\alpha}(1 + \lambda^{\alpha})
\la{sigalpha}
\ee
Also, integrating in time (\ref{tauode}) and measuring in $C^{\alpha}(\Rr^d)$ we obtain
\be
\|\tau(t_1)-\tau(t_2)\|_{C^{\alpha}(\Rr^d)}\le |t_1-t_2|D_{\alpha}
\la{taulipa}
\ee
and consequently
\be
\|\sigma(t_1)-\sigma(t_2)\|_{C^{\alpha}(\Rr^d)}\le 
|t_1-t_2|\tilde{D}_{\alpha}
\la{siglipa}
\ee
Similarly, integrating in time (\ref{tauode}) and measuring in $L^p(\Rr^d)$ we obtain
\be
\|\tau(t_1)-\tau(t_2)\|_{L^p(\Rr^d)}\le |t_1-t_2|C_p
\la{taulipp}
\ee
and similarly
\be
\|\sigma(t_1)-\sigma(t_2)\|_{L^p(\Rr^d)}\le |t_1-t_2|\tilde{C}_p
\la{siglipp}
\ee
whic we can bound using (\ref{sigp}). So we proved 
\beg{prop} \la{utosig} Let $u\in L^{\infty}(0,T; C^{1+\alpha}(\Rr^d))$ and let $\sigma(0)\in C^{\alpha, p}$ for some $\alpha\in (0,1)$ and $1<p<\infty$. Then the solution of the linear equation (\ref{sigeq}) with initial datum $\sigma(0)$ belongs to $Lip(0,T; C^{\alpha,p})$ and obeys the bounds (\ref{sigfty}), (\ref{sigp}), (\ref{sigalpha}), (\ref{siglipa}) and (\ref{siglipp}) above. Its Lagrangian counterpart $\tau$, obeys  (\ref{taufty}), (\ref{taup}), (\ref{taualpha}), (\ref{taulipa}) and (\ref{taulipp}). 
\end{prop} 
We do not use this proposition in the sequel.

Short time existence of solutions of (\ref{stokes}, \ref{sigeq}) can be  proved in the same manner as short time existence of solutions to 3D incompressible Euler equations. The stresses are Lipschitz continuous with values in $C^{\alpha, p}$ spaces, $\sigma\in Lip(0,T; C^{\alpha}(\Rr^d)\cap L^p(\Rr^d))$, and the Eulerian velocities are bounded $u\in L^{\infty}(0,T; C^{1+\alpha}(\Rr^d)\cap W^{1,p}(\Rr^d))$, for any $\alpha\in (0,1)$ and any $1<p<\infty$. Note that we do not require $p>d$. In fact, the bounds in the previous section can be used to prove a local existence theorem. In this section we investigate properties of linearizations along families of Lagrangian paths and prove existence and uniqueness results.

We take a uniformly bounded family of paths, depending in a differentiable manner of a parameter $\epsilon$,  $X_{\epsilon}\in Lip(0,T; C^{1+\alpha}(\Rr^d))$ and a uniformly bounded family depending in a differentiable manner of $\epsilon$, $\tau_{\epsilon}\in Lip(0,T; C^{\alpha, p})$ with initial data $\sigma_{\epsilon}(0)$. We assume that $X_{\epsilon}-{\mathbb I}$ is bounded in $Lip(0,T; C^{1+\alpha}(\Rr^d)\cap W^{1,p}(\Rr^d))$. We consider uniformly bounded, $\epsilon$-differentiable family of initial data, $u_{\epsilon}(0)\in C^{1+\alpha}(\Rr^d)\cap W^{1,p}(\Rr^d)$.

Measuring $\mathcal U'$ given by (\ref{xprimeeq}) in $C^{\alpha, p}$  using (\ref{xprimeeul}) and the bounds (\ref{usbound}, (\ref{ucommbound}) we obtain
\be
\ba
\left \|\mathcal U'\right \|_{L^{\infty}(0,T; C^{\alpha, p})} \\
\le \epsilon(T)\left\{\|X'\|_{C^{\beta}(0,T; C^{\alpha, p})} + \|X'\|_{L^{\infty}(0,T; C^{1+\alpha, p})} + \|\tau'\|_{L^{\infty}(0,T; C^{\alpha, p})}\right\} + C\|u'(0)\|_{C^{\alpha, p}}
\ea
\la{xprimebound}
\ee
Let us denote the norms in the right hand side of (\ref{xprimebound}) by $N(T)$:
\be
N(T)= \sup_{\epsilon}\left\{\|X'_{\epsilon}\|_{C^{\beta}(0,T; C^{\alpha, p})} + 
\|X'_{\epsilon}\|_{L^{\infty}(0,T; C^{1+\alpha, p})} + \|\tau'_{\epsilon}\|_{L^{\infty}(0,T; C^{\alpha, p})} + \|u'_{\epsilon}(0)\|_{C^{\alpha, p}}\right\}
\la{nt}
\ee
We make the convention that $C(T)$ denotes a constant that depends continuously, nondecreasingly, and explicitly on $T$, and we will use $\epsilon(T)$ for constants that vanish at least like $\max\{T^{1-\beta}, T^{\fr{1}{2}}\}$ at $T=0$. Time independent constants are written as $C$. The constants $C(T)$ , $\epsilon(T)$ depend on the assumed uniform bounds on the families $X_{\epsilon}, \tau_{\epsilon}$.
Let us introduce
\be
{\mathcal X}'(a,t) = \int_0^t{\mathcal U}'(a,s)ds.
\la{calxprime}
\ee

Integrating in time, we have from (\ref{xprimebound}) 
\be
\left \|{\mathcal X}'\right \|_{L^{\infty}(0,T; C^{\alpha, p})}\le \epsilon(T)N(T)
\la{linfty}
\ee

In order to close the estimates we have to consider a stronger path norm:
\be
M(T) =
\sup_{\epsilon}\left\{ \|X'_{\epsilon}\|_{C^{\beta}(0,T; C^{1+\alpha, p})} + 
\|\tau'_{\epsilon}\|_{C^{\beta}(0,T; C^{\alpha, p})} + \|u'_{\epsilon}(0)\|_{C^{1+\alpha, p}}\right\}.
\la{mt}
\ee
Clearly, $N(T)\le M(T)$. Let us bound $g'$ given in (\ref{gprime}), using (\ref{gprimeeul}), Theorem \ref{gbound} and Theorem \ref{commthm}. The term 
$({\mathbb L}\na\na u_0)\eta$ is bounded using
\[
\|{\mathbb L}\na\na u_0\|_{C^{\alpha, p}}\le Ct^{-\fr{1}{2}}\|u(0)\|_{C^{1+\alpha, p}}
\]
combined with the bound
\[
\|\eta(t)\|_{C^{\alpha}(\Rr^d)}\le t^{\beta}M(T)
\]
which follows from the definition (\ref{mt}) and the fact that $\eta(0)=0$.
We use here $\beta>\fr{1}{2}$. We obtain
\be
\|g'(t)\|_{C^{\alpha, p}}\le C(T)M(T).
\la{gprimebound}
\ee
Now we use (\ref{naxiprimeeq}) and (\ref{gprimebound}) to bound
\be
\|\na {\mathcal U}'\|_{L^{\infty}(0,T; C^{\alpha, p})} \le C(T)M(T).
\la{nauprimebound}
\ee
Consequently, using (\ref{nauprimebound}) in (\ref{calxprime}) we have
\be
\|\na {\mathcal X}'\|_{L^{\infty}(0,T; C^{\alpha, p})}\le TCM(T)=\epsilon(T)M(T)
\la{xprimec1alpha}
\ee
and also
\be
\| {\mathcal X}'\|_{C^{\beta}(0,T; C^{1+\alpha, p})}\le T^{1-\beta}C(T)M(T) = \epsilon(T)M(T).
\la{xprimecbeta}
\ee
We used here that an $O(1)$ bound on the time derivative gives, for short time
an $O(t^{1-\beta})$ bound on the $C^{\beta}$ norm in time (if the initial data vanishes).

We  turn to (\ref{tauprimeeq}) and bound
using (\ref{gprimebound})
\be
\|{\mathcal T}'\|_{L^{\infty}(0,T; C^{\alpha, p})}\le C(T)M(T).
\la{tauprimebound}
\ee
We define now
\be
\pi(a,t) = \int_0^t {\mathcal T}'(a,s)ds + \sigma'(0)
\la{pi}
\ee
and deduce from (\ref{tauprimebound}) 
\be
\|\pi\|_{C^{\beta}(0,T; C^{\alpha, p})}\le \epsilon(T)M(T) + \|\sigma'(0)\|_{C^{\alpha, p}}
\la{tauprimecbeta}
\ee
Summarizing what we have obtained in (\ref{xprimecbeta}) and (\ref{tauprimecbeta})
\be
\|{\mathcal X}'\|_{C^{\beta}(0,T; C^{1+\alpha, p})} + \|\pi\|_{C^{\beta}(0,T; C^{\alpha, p})}\le \epsilon(T)M(T) + \|\sigma'(0)\|_{C^{\alpha, p}}
\la{sofar}
\ee
This is the main inequality of this section. It will be used in several situations. First, let us consider the map $\mathcal S$,
\be
(X,\tau) \mapsto  {\mathcal S}(X,\tau) = (X_{new}, \tau_{new})
\la{cals}
\ee
defined by
\be
X_{new}(a,t) = a + \int_0^t{\mathcal U}(X(s), \tau (s))ds
\la{xnew}
\ee
and
\be
\tau_{new}(a,t) = \sigma_0(a) + \int_0^t{\mathcal T}(X(s),\tau(s))ds 
\la{taunew}
\ee
where ${\mathcal U}, \; {\mathcal T}$ are given in (\ref{lagsys}).
\beg{thm}\la{laglip} Let $0<\alpha<1$, $\fr{1}{2}<\beta<1$, $1<p<\infty$ and let $u_0\in C^{1+\alpha, p}$ and $\sigma_0\in C^{\alpha, p}$ be fixed. There exists $T$ sufficienytly small such that the map ${\mathcal S}$ maps the set
\be
\ba
{\mathcal I}\subset  {\mathcal P}_1 = Lip(0,T; C^{1+\alpha, p})\times Lip(0,T; C^{\alpha,p})\\
{\mathcal I} = \{(X,\tau)\left |\right.\; \left \|(X-{\mathbb I},\tau)\right \|_{{\mathcal P}_1}\le \Gamma, \; \fr{1}{2}\le |\na_a X(a,t)| \le \fr{3}{2}\}
\ea
\la{cali}
\ee
to itself, ${\mathcal S}: {\mathcal I}\to {\mathcal I}$. Furthermore, the map is a contraction in the space
\be
{\mathcal P} =C^{\beta}(0,T; C^{1+\alpha, p})\times C^{\beta}(0,T; C^{\alpha,p})
\la{calp}
\ee
i.e.,
\be
\left \|{\mathcal S}(X_1,\tau_1) -{\mathcal S}(X_2,\tau_2)\right \|_{\mathcal P}\le \fr{1}{2}\left \|(X_1-X_2,\tau_1-\tau_2)\right \|_{\mathcal P}.
\la{laglipbound}
\ee
for $(X_1,\tau_1)\in {\mathcal I}$, $(X_2, \tau_2)\in \mathcal I$.
\end{thm}
\noindent{\bf Proof.} The fact that ${\mathcal S}:{\mathcal I}\to {\mathcal I}$
follows from the bounds in Theorem \ref{uboun} and Theorem \ref{gbound} by choosing $T$ small enough, and $\Gamma$ twice the size of the initial data $(u_0, \sigma_0)$ in $C^{1+\alpha, p} \times C^{\alpha, p}$. The contractivity is proved by forming the families $X_{\epsilon} = (2-\epsilon)X_1 +(\epsilon -1)X_2$, $\tau_{\epsilon} = (2-\epsilon)\tau_1 + (\epsilon -1)\tau_2$ for $\epsilon\in[1,2]$. We note that ${\mathcal I}$ is convex, and that $u'_0 =0$, because $u_0$ does not depend on epsilon, and $\sigma'_0 = 0$ because $\sigma_0$ does not depend on epsilon. 
\be
{\mathcal S}(X_1,\tau_1)-{\mathcal S}(X_2,\tau_2) = \left(\int_1^2{\mathcal X}'_{\epsilon}d\epsilon, \int_1^2\pi_{\epsilon}d\epsilon\right)
\la{diff}
\ee
where ${\mathcal X}'$ and $\pi$ are obtained via (\ref{calxprime}) and (\ref{pi}) using the families $X_{\epsilon},\tau_{\epsilon}$. Applying (\ref{sofar}) and choosing $T$ small enough proves the contractivity. The local existence and uniqueness theorem is:

\beg{thm}\la{lagexist} Let $0<\alpha<1$, $1<p<\infty$, $\fr{1}{2}<\beta<1$, and let $u_0\in C^{1+\alpha, p}$ be divergence-free, and $\sigma_0\in C^{\alpha, p}$ be given.

\noindent (A) There exists $T>0$ and a solution $(u,\sigma)$ of (\ref{stokes}), (\ref{sigeq})with $u\in L^{\infty}(0,T; C^{1+\alpha, p})$ and $\sigma\in Lip(0,T; C^{\alpha,p})$.

\noindent (B) Two solutions $u_j\in L^{\infty}(0,T; C^{1+\alpha,p})$ and $\sigma_j\in Lip(0,T; C^{\alpha,p})$, $j=1,2$ obey the strong Lipschitz bound
\be
\ba
\|X_2-X_1\|_{C^{\beta}(0,T; C^{1+\alpha, p})} + \|\tau_2-\tau_1\|_{C^{\beta}(0,T; C^{\alpha, p})} \le \\
C(T)\left \{\|u_2(0)-u_1(0)\|_{C^{1+\alpha, p}} +
 \|\sigma_2(0)-\sigma_1(0)\|_{\alpha,p}\right\}
\ea,
\la{strong}
\ee
for their Lagrangian counterparts. The path time derivatives also obey Lipschitz bounds
\be
\|\pa_t X_2-\pa_t X_1\|_{L^{\infty}(0,T; C^{1+\alpha,p})}\le C(T)\left \{\|u_2(0)-u_1(0)\|_{C^{1+\alpha, p}} + \|\sigma_2(0)-\sigma_1(0)\|_{\alpha,p}\right\}
\la{patx}
\ee
In particular, two such solutions with same initial data must coincide.
\end{thm}

\noindent{\bf{Proof}}. Part (A), the existence, follows because a fixed point of ${\mathcal S}$ provides a solution in ${\mathcal I}$. The initial velocity being divergence-free and the equation (\ref{u}) guarantee incompressibility. Part
(B) is proved forming the family $(X_{\epsilon}, \tau_{\epsilon})$ as in the proof of Theorem \ref{laglip} above, with $(X_1, \tau_1)$ being the Lagrangian solution associated to the solution $(u_1,\sigma_1)$ and  with $(X_2,\tau_2)$ being the Lagrangian solution associated to the solution $(u_2,\sigma_2)$.
\[
X_{\epsilon} = (2-\epsilon)X_1 +(\epsilon -1)X_2, \quad \tau_{\epsilon} = (2-\epsilon)\tau_1 + (\epsilon -1)\tau_2
\]
for $\epsilon\in[1,2]$. Note that $X'_{\epsilon} = X_2-X_1$, $\tau'_{\epsilon} = \tau_2-\tau_1$ but also, because these are solutions,
\be
X_2-X_1 = \int_1^2 {\mathcal X}'_{\epsilon}d\epsilon 
\la{difx}
\ee
and
\be
\tau_2-\tau_1 =\int_1^2\pi_{\epsilon}d\epsilon + \sigma_2(0)-\sigma_1(0) 
\la{diftau}
\ee
Also, $u'_0 = u_2(0)-u_1(0)$, $\sigma'_0 = \sigma_2(0)-\sigma_1(0)$. So we have
\be
\left\{
\ba
\int_1^2{\mathcal X}'_{\epsilon}d\epsilon = X_2-X_1= X',\\
\int_1^2\pi_{\epsilon}d\epsilon = \tau_2-\tau_1-\sigma_2(0)+ \sigma_1(0) = \tau'-\sigma_2(0)+ \sigma_1(0).
\ea
\right.
\la{together}
\ee
Integrating in epsilon (\ref{sofar}) we obtain
\be
\ba
\|X_2-X_1\|_{C^{\beta}(0,T; C^{1+\alpha, p})} + \|\tau_2-\tau_1-\sigma_2(0)+\sigma_1(0)\|_{C^{\beta}(0,T; C^{\alpha, p})} \le \\
\epsilon(T)\left \{\|X_2-X_1\|_{C^{\beta}(0,T; C^{1+\alpha, p})} + \|\tau_2-\tau_1\|_{C^{\beta}(0,T; C^{\alpha, p})} + \|u_2(0)-u_1(0)\|_{C^{1+\alpha, p}}\right\} \\
+ \|\sigma_2(0)-\sigma_1(0)\|_{\alpha,p}
\ea
\la{semifin}
\ee 
Taking $\epsilon(T)\le \fr{1}{2}$ we obtain the strong Lipschitz bound (\ref{strong}). Note also that this implies 
\be
M(T)\le C(T)\left \{\|u_2(0)-u_1(0)\|_{C^{1+\alpha, p}} + \|\sigma_2(0)-\sigma_1(0)\|_{\alpha,p}\right\}
\la{temp}
\ee
and therefore, going back to (\ref{xprimebound}) and (\ref{nauprimebound}), and using
\[
\pa_tX_2-\pa_t X_1 = \int_1^2 {\mathcal U}'_{\epsilon}d\epsilon
\]
we obtain (\ref{patx}).

\section{Coupling to Navier-Stokes}
The Navier-Stokes equations are nonlinear
\be
\pa_t u -\Delta u + \na p =\div (\sigma - u\otimes u), \quad \div u = 0.
\la{nseq}
\ee
In order to prove uniqueness of solutions of the system formed by (\ref{nseq}) coupled to (\ref{sigeq}) we still have to work in a class of velocities that are at least Lipschitz continuous. This is a vastly subcritical situation for Navier-Stokes equations, so we treat the inertial stress $\div(u\otimes u)$ perturbatively. We write
\be
u  = {\mathbb L}(u_0) + {\mathbb U}(\sigma) - {\mathbb U}(u\otimes u)
\la{unse}
\ee
and
\be
\na u = {\mathbb L}(\na u_0) + {\mathbb G}(\sigma) - {\mathbb U}(\na(u\otimes u)).
\la{nauns}
\ee
We wrote ${\mathbb U}(\na(u\otimes u))$ above instead of the equivalent ${\mathbb G}(u\otimes u)$, in order to take advantage of the fact that in our framework $\na u\in C^{\alpha, p}$.
We introduce again the Lagrangian variables $X$ and $\tau$, but we keep a separate tab for the Lagrangian velocity
\be
v= u\circ X.
\la{vns}
\ee
We set up the map
\be
{\mathcal V}(X,\tau, v) = {\mathbb L}(u_0)\circ X + {\mathbb U}(\tau\circ X^{-1})\circ X - {\mathbb U}((v\otimes v)\circ X^{-1})\circ X,
\la{mathcalvnse}
\ee
we denote again
\be
g= {\mathbb L}(\na u_0)\circ X + {\mathbb G}(\tau\circ X^{-1})\circ X - {\mathbb U}(\na((v\otimes v)\circ X^{-1}))\circ X
\la{gnse}
\ee
and we maintain
\be
{\mathcal T}(X,\tau, v) = F(g,\tau).
\la{tnse}
\ee
The system is now solved in Lagrangian coordinates:
\be
(X,\tau, v)\mapsto {\mathcal S}(X,\tau, v) = (X^{new}, \tau^{new}, v^{new})
\la{mathcalsns}
\ee
where
\be
\left\{
\ba
X^{new}(t) = a + \int_0^t{\mathcal V}(X(s),\tau(s), v(s))ds\\
\tau^{new} (t) = \sigma_0 + \int_0^t{\mathcal T}(X(s), \tau(s), v(s))ds\\
v^{new}(t) = {\mathcal V}(X,\tau, v)
\ea
\right.
\la{nsxtauveq}
\ee
Note that
\be
{\mathcal V}(X,\tau, v) = {\mathcal U}(X,\tau) - {\mathbb U}((v\otimes v)\circ X^{-1})\circ X
\la{VU}
\ee
where ${\mathcal U}$ is the same as the one given in (\ref{lagsys}). Note also that
\be
v^{new} = \fr{dX^{new}}{dt}
\la{vxt}
\ee
and therefore the relation  $v =\fr{dX}{dt}$ is maintained in an iteration.
Because the last equation of (\ref{nsxtauveq}) is not integrated in time, we measure $v$ in $L^{\infty}(0,T; C^{1+\alpha}(\Rr^d)\cap W^{1,p}(\Rr^d))$. 
We take again a family $X_{\epsilon}, \tau_{\epsilon}, v_{\epsilon}$, denote
\be
\fr{d}{d\epsilon}v_{\epsilon} = v'_{\epsilon}
\la{vprime}
\ee
and keep the rest of the notation $X'_{\epsilon}, \eta_{\epsilon}, \tau'_{\epsilon},\delta_{\epsilon}$ the same.
Note that
\be
v'_{\epsilon} = \fr{dX'_{\epsilon}}{dt}
\la{vprimexprimet}
\ee
We differentiate the nonlinearities ${\mathcal V}$ and ${\mathcal T}$ with respect to epsilon. After composition with $X_{\epsilon}^{-1}$ and dropping epsilon for ease of notation, we have 
\be
\left\{
\ba
{\mathcal V}'\circ X^{-1} = {\mathcal U}'\circ X^{-1} - [\eta\cdot\na, {\mathbb U}](u\otimes u) - {\mathbb U}(((v\otimes v')+ (v'\otimes v))\circ X^{-1})\\
{\mathcal T}' = D_1F(g,\tau)g'+D_2F(g,\tau)\tau'\\
g'\circ X^{-1} = g'_{old}\circ X^{-1} - [\eta\cdot\na, {\mathbb U}](\na(u\otimes u))-{\mathbb{U}}(\na(((v'\otimes v) + (v\otimes v'))\circ X^{-1}))
\ea
\right.
\la{lagssyprimens}
\ee 
where ${\mathcal U}'$ and $g'_{old}$ are the same as in (\ref{linlageq}) and $u= v\circ X^{-1}$. 
We verify that ${\mathcal S}$ maps the set
\be
\ba
{\mathcal I} \subset {\mathcal P}_1 = Lip(0,T; C^{1+\alpha, p})\times Lip(0,T; C^{\alpha,p})\times L^{\infty}(0,T; C^{1+\alpha, p})\\
{\mathcal I} = \{(X,\tau, v)\left |\right. \; \|(X-{\mathbb I}, \tau , v)_{{\mathcal P}_1} \le \Gamma, \fr{1}{2}\le |\na_a X(a,t)|\le \fr{3}{2}, v =\fr{dX}{dt}\}
\ea
\la{mathcali}
\ee
to itself for $\Gamma$ larger than the size of the initial data $\sigma_0$, $u_0$ and small enough $T$. In order to check this, we use the bounds used in the previous section. In addition, in view of (\ref{usbound}) we see that
\be
\|{\mathbb U}((v\otimes v)\circ X^{-1})\circ X\|_{L^{\infty}(0,T; C^{1+\alpha,p})}\le \sqrt{T}C\|v\|_{L^{\infty}(0,T; C^{1+\alpha,p})}^2
\la{uvbound}
\ee
for invertible $X$ satisfying the constraint $\fr{1}{2}\le |\na_a X(a,t)|\le 
\fr{3}{2}$. 

The new terms introduced in ${\mathcal V}'$ and in $g'$ are bounded using (\ref{usbound}) and Theorem \ref{ucomm}. Using (\ref{xprimebound}) we get
\be
\|{\mathcal V}'\|_{L^{\infty}(0,T; C^{\alpha, p})} \le
\epsilon(T)\left\{\|X'\|_{C^{\beta}(0,T; C^{1+\alpha,p})} + \|\tau'\|_{C^{\beta}(0,T; C^{\alpha,p})} + \|v'\|_{L^{\infty}(0,T; C^{\alpha,p})}\right\} + C\|u'_0\|_{C^{\alpha,p}}.
\la{vprimebound}
\ee

In order to bound $g'\circ X^{-1}$ we use (\ref{gprimebound}) to bound the term $g'_{old}\circ X^{-1}$, then (\ref{ucommbound}) together with the uniform bound on $\na(u\otimes u)$ in order to bound  $[\eta\cdot\na, {\mathbb U}](\na (u\otimes u))$, and (\ref{usbound}) for the last term in $g'\circ X^{-1}$.
We obtain
\be
\|g'\|_{L^{\infty}(0,T; C^{\alpha,p})} \le C(T)M(T) + \epsilon(T)\|v'\|_{L^{\infty}(0,T; C^{1+\alpha,p})}
\la{gprimeboundns}
\ee
where $M(T)$ is given by (\ref{mt}).
The relation 
\be
\na_a {\mathcal V}(a,t) = g(a,t)(\na_a X(a,t))
\la{naav}
\ee
is directly verified. Differentiating in epsilon, we obtain
\be
(\na_a {\mathcal V}')(a,t) = g'(a,t)(\na_a X)(a,t) + g(a,t)(\na_a X')(a,t)
\la{naavprime}
\ee
and then we obtain using (\ref{vprimebound}) and (\ref{gprimeboundns})
\be
\|{\mathcal V}'\|_{L^{\infty}(0,T; C^{1+\alpha, p})}\le 
\epsilon(T)\|v'\|_{L^{\infty}(0,T; C^{1+\alpha,p})} + C(T)M(T).
\la{vprimegradbound}
\ee
We recall that $\pi$ is given by (\ref{pi}) and introduce
\be
{\mathcal X}'(a,t) = \int_0^t{\mathcal V}'(a,s)ds.
\la{mathcalxprimens}
\ee
We have therefore
\be
\|{\mathcal X}'\|_{C^{\beta}(0,T; C^{1+\alpha, p})} + 
\|\pi\|_{C^{\beta}(0,T; C^{\alpha, p})}\le \epsilon(T)M_1(T) + \|\sigma'(0)\|_{C^{\alpha, p}}
\la{sofarns}
\ee
where
\be
M_1(T) = \sup_{\epsilon}\left\{ \|X'_{\epsilon}\|_{C^{\beta}(0,T; C^{1+\alpha, p})} + \|\tau'_{\epsilon}\|_{C^{\beta}(0,T; C^{\alpha, p})} + 
\|v'_{\epsilon}\|_{L^{\infty}(0,T; C^{1+\alpha,p})} + \|u'_{\epsilon}(0)\|_{C^{1+\alpha, p}}\right\}
\la{m1t}
\ee
The inequalities (\ref{vprimegradbound}) and (\ref{sofarns}) are the main inequalities of this section. They can be used to prove local existence and Lipschitz dependence on initial data in Lagrangian variables.
\beg{thm}\la{contrns} Let $0<\alpha<1$, $1<p<\infty$, $\fr{1}{2}<\beta<1$, let $u_0\in C^{1+\alpha, p}$ be divergence-free and let $\sigma_0\in C^{\alpha, p}$. Then there exists $\Gamma>0$, $T>0$, $\delta>0$ such that the map $\mathcal S$ defined in 
(\ref{nsxtauveq}) maps ${\mathcal I}$ defined in (\ref{mathcali}) to itself. Moreover
\be
\|{\mathcal S}(X_1,\tau_1, v_1)-{\mathcal S}(X_2,\tau_2, v_2)\|_{\mathcal P}\le \fr{1}{2}\|(X_1,\tau_1,v_1)-(X_2,\tau_2, v_2)\|_{\mathcal P}
\la{lipns}
\ee
where
\be
%\ba
\|(X_1,\tau_1,v_1)-(X_2,\tau_2, v_2)\|_{\mathcal P} =
\|X_1-X_2\|_{C^{\beta}(0,T; C^{1+\alpha,p})} +
\|\tau_1-\tau_2\|_{C^{\beta}(0,T; C^{\alpha,p})} + \delta \|v_1-v_2\|_{L^{\infty}(0,T; C^{1+\alpha,p})}
%\ea
\la{pnorm}
\ee

\end{thm}
\noindent{\bf{Proof}}. We already gave sufficient grounds to verify the fact that ${\mathcal S}$ maps $\mathcal I$ to itself for appropriate $\Gamma$ and $T$. In order to verify the contraction property, given a pair $(X_j, \tau_j, v_j=\fr{dX_j}{dt})$, $j=1,2$, we form the family $X_{\epsilon} = (2-\epsilon)X_1 + (\epsilon-1)X_2$, $\tau_{\epsilon} = (2-\epsilon)\tau_1 + (\epsilon-1)\tau_2$ and $v_{\epsilon} = (2-\epsilon)v_1 +(\epsilon-1)v_2$ and use the fact that
\[
{\mathcal S}(X_2, \tau_2, v_2)- {\mathcal S}(X_1, \tau_1, v_1) = \left(\int_1^2{\mathcal X'}_{\epsilon}d\epsilon, \int_1^2\pi_{\epsilon}d\epsilon, \int_{1}^2{\mathcal V'}_{\epsilon}d\epsilon\right)
\]
Now $u'_0=0$ and $\sigma'(0) = 0$ because the initial data are fixed and therefore
\[
M(T) = \|X_2-X_1\|_{C^{\beta}(0,T; C^{1+\alpha, p})} + \|\tau_1-\tau_2\|_{C^{\beta}(0,T; C^{\alpha,p})}
\]
and
\[
M_1(T) = M(T) + \|v_1-v_2\|_{L^{\infty}(0,T; C^{1+\alpha,p})}. 
\]
From (\ref{vprimegradbound}) we deduce
\be
\|v_2^{new}-v_1^{new}\|_{L^{\infty}(0,T; C^{1+\alpha, p})}\le \epsilon(T)\|v_1-v_2\|_{L^{\infty}(0,T; C^{1+\alpha,p})} + C(T)M(T)
\la{bad}
\ee
and from (\ref{sofarns}) we deduce
\be
\|X^{new}_1-X^{new}_2\|_{C^{\beta}(0,T; C^{1+\alpha,p})} +
\|\tau^{new}_1-\tau^{new}_2\|_{C^{\beta}(0,T; C^{\alpha,p})}\le \epsilon(T)M_1(T)
\la{good}
\ee
Let us fix $C(T)=K$ in (\ref{vprimegradbound}) and let $T$ be small enough so that $\epsilon(T)\le{\fr{1}{4}}$ and $\epsilon(T)\le \fr{1}{2(1+4K)}$. We choose $\delta =\fr{1}{4K}$, multiply (\ref{bad}) by $\delta$ and add to (\ref{good}): we obtain (\ref{lipns}).

\beg{thm}\la{exns} Let $0<\alpha<1$, $1<p<\infty$, $\fr{1}{2}<\beta<1$, let $u_0\in C^{1\
+\alpha, p}$ be divergence-free and let $\sigma_0\in C^{\alpha, p}$. 

\noindent (A) There exists $T>0$ and a solution $(u,\sigma)$ of the system (\ref{nseq}), (\ref{sigeq}) with $u\in L^{\infty}(0,T; C^{1+\alpha, p})$ and with $\sigma\in Lip(0,T; C^{\alpha,p})$.

\noindent (B) Two solutions $u_j\in L^{\infty}(0,T; C^{1+\alpha,p})$ and $\sigma_j\in Lip(0,T; C^{\alpha,p})$, $j=1,2$ obey the strong Lipschitz bound
\be
\ba
\|X_2-X_1\|_{C^{\beta}(0,T; C^{1+\alpha, p})} + \|\tau_2-\tau_1\|_{C^{\beta}(0,T; C^{\alpha, p})} +\|\pa_t X_2-\pa_t X_1\|_{L^{\infty}(0,T; C^{1+\alpha, p})}\le \\
C(T)\left \{\|u_2(0)-u_1(0)\|_{C^{1+\alpha, p}} + \|\sigma_2(0)-\sigma_1(0)\|_{\alpha,p}\right\}
\ea
\la{strongns}
\ee
for their Lagrangian counterparts. In particular, two such solutions with the same initial data must coincide.

\end{thm}
\noindent{\bf Proof.} The proof is very similar to the proof of Theorem \ref{lagexist}. Using (\ref{together}) and integrating in epsilon in (\ref{sofarns}) we obtain
\be
\ba
\|X_2-X_1\|_{C^{\beta}(0,T; C^{1+\alpha, p})} + \|\tau_2-\tau_1-\sigma_2(0)+\sigma_1(0)\|_{C^{\beta}(0,T; C^{\alpha, p})} \le \\
\epsilon(T)\left \{\|X_2-X_1\|_{C^{\beta}(0,T; C^{1+\alpha, p})} + \|\tau_2-\tau_1\|_{C^{\beta}(0,T; C^{\alpha, p})}  +\|v_2-v_1\|_{L^{\infty}(0,T; C^{1+\alpha, p})}\right\}\\ + \epsilon(T)\|u_2(0)-u_1(0)\|_{C^{1+\alpha, p}}
+ \|\sigma_2(0)-\sigma_1(0)\|_{\alpha,p}. 
\ea
\la{semifins}
\ee 
On the other hand, using
\[
v_2-v_1 = \int_1^2 {\mathcal V}'_{\epsilon}d\epsilon
\]
and integrating with respect to epsilon in (\ref{vprimegradbound}) we obtain
\be
\|v_2-v_1\|_{L^{\infty}(0,T; C^{1+\alpha, p})}\le \epsilon(T)\|v_2-v_1\|_{L^{\infty}(0,T; C^{1+\alpha, p})} + C(T)M(T)
\la{semifinv}
\ee
Fixing $C(T)= K$ and choosing again $T$ small enough such that $\epsilon(T)\le \fr{1}{4}$ and $\epsilon(T)\le \fr{1}{2(1+4K)}$ we obtain that
\be
\ba
\|X_2-X_1\|_{C^{\beta}(0,T; C^{1+\alpha, p})} + \|\tau_2-\tau_1\|_{C^{\beta}(0,T; C^{\alpha, p})} + \fr{1}{4K}\|v_2-v_1\|_{L^{\infty}(0,T; C^{1+\alpha, p})}\\
\le C(T)\left \{\|u_2(0)-u_1(0)\|_{C^{1+\alpha, p}} + \|\sigma_2(0)-\sigma_1(0)\|_{\alpha,p}\right\}
\ea
\la{fins}
\ee
and therfore (\ref{strongns}) is proved.

\subsection*{Acknowledgment}Research partially supported by NSF-DMS grants 1209394 and 1265132.


\begin{thebibliography}{99}
\bibitem{csun}P. Constantin, W. Sun, Remarks on Oldroyd-B and related complex fluid models, CMS, {\bf{10}} No. 1, (2012), 33-73.
\bibitem{fanghuachun} F. Lin, C. Liu, P. Zhang, On hydrodynamics of viscoelastic fluids, Comm. Pure  Appl. Math {\bf{58}} (2005), 1437-1471.

\bibitem{mb} A. Majda, A. Bertozzi, {\em{Vorticity and Incompressible Flow}}, Cambridge Texts in Appl. Math, CUP, Cambridge, 2002.

\bibitem{renardy} M. Renardy, An existence theorem for model equations resulting from kinetic theories of polymer solutions, SIAM J. Math. Anal., {\bf{22}} (1991), 3131-327.

\bibitem{saut} C. Guillop\'{e}, J.-C. Saut, Existence results for the flow of viscoelastic fluids with a differential constitutive law, Nonlinear Anal., {\bf{15}} (1990), 849-869.

\end{thebibliography}
\end{document}